\documentclass{amsart}

\usepackage{graphicx}
\usepackage{pst-all}      
\usepackage{url}

\newtheorem{theorem}{Theorem}

\newtheorem{question}{Question}

\begin{document}

\title{Generating irreducible triangulations of surfaces}


\author{Thom Sulanke}
\address{Department of Physics, Indiana University, Bloomington, 
	Indiana 47405}
\email{tsulanke@indiana.edu}

\date{\today}

\begin{abstract}
Starting with the irreducible triangulations of a fixed
surface and splitting vertices, all the triangulations of the 
surface up to a given number of vertices can be generated.
The irreducible triangulations have previously been determined for 
the surfaces $S_0$, $S_1$, $N_1$,and $N_2$.
An algorithm is presented for generating the irreducible 
triangulations of a fixed surface using triangulations of other surfaces.
This algorithm has been implemented as a computer program
which terminates for $S_1$, $S_2$, $N_1$, $N_2$, $N_3$, and $N_4$. 
Thus the complete sets irreducible triangulations are now also known for 
$S_2$, $N_3$, and $N_4$, with respective cardinalities 396784, 9708, and 
6297982.
\end{abstract}

\maketitle

\section{Introduction}
\label{intro}

We can transform a triangulation of a surface into another triangulation of the
same surface by cutting along two edges which have a common end point and 
inserting two new triangular faces.
By repeatedly applying this operation of vertex splitting to one of the small 
(irreducible) triangulations of a surface we can obtain any given 
triangulation of that surface.

We can transform a triangulation of a surface into a triangulation of a 
different surface in much the same way we transform one surface into a
different surface by adding a handle, a crosshandle, or a crosscap.
To add a handle to a triangulation we remove two faces from the triangulation
and identify the two new boundaries.
If these two faces are far enough apart then no loops and no multiple edges are
produced and we have a triangulation of the original surface with a 
handle added.
To add a crosscap to a triangulation we remove a vertex of degree six along 
with the edges and faces incident to it.
This produces a hexagonal hole.  
We produce a crosscap by identifying the vertices of the three pairs of 
opposite vertices of this hole.
If the vertices in these pairs are far enough apart then no loops and no 
multiple edges are produced and we have a triangulation of the 
original surface with a crosscap added.

We show in Section~\ref{genirrtri} that every irreducible triangulation of a 
surface other than the sphere can be obtained by adding a handle, crosshandle,
or crosscap to another triangulation.
Thus for any triangulation of a surface there is a sequence of vertex 
splittings, handle additions, crosshandle additions, and crosscap additions 
which transform the tetrahedron into the given triangulation.
For a fixed surface there are a finite number of irreducible triangulations.
Therefore, the irreducible triangulations of a fixed surface can be produced
in a finite number of steps.
We present an algorithm which not only produces all the irreducible 
triangulations of $S_1$, $S_2$, $N_1$, $N_2$, $N_3$, and $N_4$, 
but also terminates after a finite number of steps for these surfaces.

In Section~\ref{verify} we discuss several procedures that were used to confirm
the results obtained by computer programs which implement the algorithms 
presented in earlier sections.

The irreducible triangulations of the sphere (Steinitz and Rademacher 
\cite{StRa}), the projective plane
(Barnette \cite{MR84f:57009}), the torus (Lawrencenko \cite{MR914777}), and the
Klein bottle (Lawrencenko and Negami \cite{MR98h:05067} and 
Sulanke \cite{math.CO/0407008}) have been determined. 
The number of irreducible triangulations for each surface is
1 for $S_0$, 2 for $N_1$, 21 for $S_1$, and 29 for $N_2$.
Using the algorithm described here the number of irreducible 
triangulations have been found to be 9708 for $N_3$, 396784 for $S_2$, and 
6297982 for $N_4$.

\section{Definitions}
\label{definitions}

A {\em triangulation\/} of a closed surface 
is a simple graph embedded in the surface
such that each face is a triangle and any two faces share at most one edge.   

In a triangulation $T$ let
$abc$ and $acd$ be two faces which have $ac$ as a common edge.
The {\em contraction\/} of $ac$ is obtained by deleting $ac$,
identifying vertices $a$ and $c$, 
removing one of the multiple edges $ab$ or $cb$,
and removing one of the multiple edges $ad$ or $cd$.
The edge $ac$ of a triangulation $T$ is {\em contractible\/} if the 
contraction of $ac$ yields another triangulation of the surface in
which $T$ is embedded. 
If the edge $ac$ is contained in a 3-cycle other than the two 
which bound the faces which share it then its contraction 
would produce multiple edges.
Thus, for a triangulation $T$, not $K_4$ embedded in the sphere, 
an edge of $T$ is not contractible
if and only if that edge is contained in at least three 3-cycles. 
A vertex is said to be {\em contractible\/} 
if it is the end of at least one contractible edge. 
A triangulation is said to be {\em irreducible\/} 
if it has no contractible edges. 

Let $abc$ be a 3-cycle in $T$ which is embedded in the surface $S$ and 
let $C$ be 
the closed curve which is the embedding of $abc$ in $S$.
Then $abc$ is said to be a {\em separating\/} 3-cycle if $S - C$ is 
disconnected and
$abc$ is said to be a {\em nonseparating\/} 3-cycle, otherwise.
If the neighborhood of $C$ in $S$ is homeomorphic to a M\"{o}bius band then
we say that $abc$ is a {\em one-sided\/} 3-cycle, 
otherwise we say that $abc$ is a {\em two-sided\/} 3-cycle.
If $abc$ bounds a face then $abc$ is a {\em facial\/} 3-cycle.
If $T$ is an irreducible triangulation of a surface other than the sphere, then
every edge of $T$ is noncontractible, and therefore on at least three 3-cycles,
two facial and the others nonfacial.

For a vertex $v$ in a triangulation we define the {\em link\/} of $v$, lk($v$),
as the cycle through the neighbors of $v$ which contains the edges of the faces
incident on $v$.

The inverse of contracting an edge is splitting a vertex.
Let $a$ be a vertex of a triangulation and let $ab$ and $ac$ be edges.
We split $a$ by cutting along the path $bac$ and inserting two new faces.
Label the neighbors of $a$ such that lk($a$) = $(v_1,v_2,\ldots,v_d)$ 
with $b=v_1$ and $c=v_s$.
The {\em splitting} of the vertex $a$ (along $bac$) is obtained by removing
$a$ along with all the edges and faces incident to $a$ and then adding 
a vertex $a_1$, 
a vertex $a_2$, 
the edge $a_1 a_2$,
the edges $\{a_1 v_i: i=1,\ldots,s\}$,
the edges $\{a_2 v_i: i=s,\ldots,d,1\}$, 
the faces $\{a_1 v_i v_{i+1}: i=1,\ldots,s-1\}$,
the faces $\{a_2 v_i v_{i+1}: i=s,\ldots,d-1\}$,
the face $a_2 v_d v_1$,
the face $a_1 a_2 v_1$, and
the face $a_1 a_2 v_s$.

We denote the orientable surface with genus $g$, the sphere with $g$ handles
attached, as $S_g$ and the 
nonorientable surface with genus $g$, the sphere with $g$ crosscaps attached,
 as $N_g$.
Define $\mathrm{eg}(S) = 2 - \chi(S)$ to be the {\em Euler genus} of the 
surface $S$.
Define the {\em Euler genus} of the triangulation $T$ 
as $\mathrm{eg}(T) = \mathrm{eg}(S)$ where $T$ is a triangulation of $S$.
For orientable surfaces $\mathrm{eg}(S_g) = 2g$ and for
nonorientable surfaces $\mathrm{eg}(N_g) = g$.


\section{Generating triangulations}
\label{gentri}

We will specify an algorithm for generating triangulations as 

\begin{itemize}

\item $\mathcal{B}$, a basic set of triangulations from which the 
triangulations of $\mathcal{T}$ are generated; 

\item $r$, a generating rule defined on $\mathcal{T}$;

\item $\overline{r}$, the inverse of the generating rule, $r$;

\item $\mathcal{T}$, the set of all triangulations generated by applying $r$ 
zero or more times to elements of $\mathcal{B}$;

\item $f$, a filter which selects some of the generated triangulations from 
$\mathcal{T}$;

\item $\mathcal{F}$, the triangulations filtered by $f$.

\end{itemize}

Thus, to obtain all the triangulations with exactly $n$ vertices of a surface 
$S$

\begin{itemize}

\item $\mathcal{B}_n$ is the set of irreducible 
triangulations of $S$ with at most $n$ vertices;

\item $r_n$ is the generating rule ``if the triangulation has less than $n$ 
vertices split a vertex'';

\item $\overline{r}_n$ is ``if the triangulation has a contractible edge then
contract a contractible edge'';

\item $\mathcal{T}_n$ is the set of 
all triangulations of $S$ with at most $n$ vertices;

\item $f_n$ is the condition ``the triangulation has exactly $n$ vertices'';

\item $\mathcal{F}_n$ is the set of all the triangulations with exactly $n$ 
vertices.

\end{itemize}

To see that $\mathcal{T}_n$ is inductively defined by $\mathcal{B}_n$ and $r_n$
consider $\overline{r}_n$, applied to elements of $\mathcal{T}_n$.
If $T \in \mathcal{T}_n - \mathcal{B}_n$ then $T$ has a contractible edge and
$\overline{r}_n$ can be applied to obtained another element of $\mathcal{T}_n$.
The number of vertices is decreased by applying $\overline{r}_n$ 
therefore, it can be 
applied only a finite number of times before an element of $\mathcal{B}_n$ is 
obtained.
Reversing these steps provides a way of obtaining $T$ from an element of 
$\mathcal{B}_n$ by applying $r_n$ sequentially.

For $T \in \mathcal{T}_n$, $r_n$ can be applied in a finite number of ways.
The number of vertices increases with each application of $r_n$ therefore, it
can be sequentially applied less than $n$ times. 
Therefore, $\mathcal{T}_n$ is finite.

Thus if we know $\mathcal{B}_n$ we have an algorithm which generates 
$\mathcal{F}_n \subset \mathcal{T}_n$ and the algorithm terminates.

This algorithm has been implemented as the author's computer program 
{\em surftri\/} \cite{surftri} using the McKay
orderly generation principle, which is described in \cite{MR1606516} and 
\cite{MR1762931}.
The program {\em surftri\/} is based extensively on the ideas, code, and data 
structures of the program {\em plantri\/} \cite{plantri}.

\section{Generating irreducible triangulations}
\label{genirrtri}

A basic step in the construction of the irreducible triangulations of 
the torus \cite{MR914777} and the
Klein bottle \cite{MR98h:05067}
is to show that any irreducible triangulation of these surfaces can be cut
along one or two 3-cycles to produce a
punctured surface of lower Euler genus.
The algorithm presented here uses a similar construction.

We examine how an irreducible triangulation can be reduced to an irreducible
triangulation of lower genus.
We then can define generating rules to reverse these steps.
Briefly we describe this reduction as follows.
Let $T$ be an irreducible triangulation of a surface other than $S_0$ with
a nonseparating 3-cycles $W$.
Cut $T$ along this 3-cycle and cap the resulting hole or holes with new
triangular faces to produce a new triangulation $T'$ of a surface with lower
Euler genus. 
Consider what can happen to a nonfacial 3-cycle $C$ in $T$ during the 
construction of $T'$.
If $C = W$ then $C$ becomes part of the cap in $T'$.
If $C$ does not ``cross'' $W$ then the edges of $C$ are also noncontractible
in $T'$.
If $C$ does ``cross'' $W$ in $T$ then $C$ becomes a path of length 3 in $T'$.
Some of the edges in the cap and some edges on these paths might be 
contractible.
The contractible edges are contracted until we obtain an irreducible 
triangulation of the surface of which $T'$ is a triangulation.

The following theorem is similar to Lemma 4 of \cite{MR84f:57009} and 
Lemma 4 of \cite{MR914777}.

\begin{theorem}
\label{transverse}
Let $T$ be an irreducible triangulation of a surface other than $S_0$,
let $v$ be a vertex of $T$.
Then there are two nonseparating 3-cycles $v v_i v_k$ and $v v_j v_l$
such that $v_i$, $v_j$, $v_k$, and $v_l$ are distinct
and one path from $v_i$ to $v_k$ in {\rm lk}($v$) contains $v_j$ and
the other path from $v_i$ to $v_k$ in {\rm lk}($v$) contains $v_l$.
\end{theorem}

\begin{proof}
Since $T$ is irreducible,
for any vertex $u$ in lk($v$) the edge $v u$ is on a nonfacial 3-cycle $v u w$.
Pick two vertices $v_i$ and $v_k$ in lk($v$) for which $v v_i v_k$ is a
nonfacial 3-cycle and the distance from $v_i$ to $v_k$ in lk($v$) is minimal.
The shorter path from $v_i$ to $v_k$ in lk($v$) must have an interior vertex
since $v v_i v_k$ is not a face.
Let the vertex $v_j$ be such an interior vertex on the shorter path from $v_i$
to $v_k$ in lk($v$).
Let $v_l$ be a vertex in lk($v$) such that $v v_j v_l$ is a nonfacial 3-cycle.
$v_l$ is not on the path from $v_i$ to $v_k$ in lk($v$) containing $v_j$ 
since the distance from $v_j$ and $v_l$ in lk($v$) is at least the distance 
from $v_i$ and $v_k$ in lk($v$).
Suppose $v v_i v_k$ separates the surface.  
Then $v_j$ and $v_l$ would be in different components but $v_j v_l$ is 
an edge.
Therefore, $v v_i v_k$ is nonseparating and, similarly, $v v_j v_l$ is also
nonseparating.
\end{proof}

Let $S_f \neq S_0$ be the surface for which we are generating irreducible 
triangulations.
Let $T$ be an irreducible triangulation of $S_f$.
By Theorem~\ref{transverse} there is a nonseparating 3-cycle $w_1 w_2 w_3$ 
in $T$.
We create a new triangulation $T'$ of a different surface $S_b$ 
using one of the two following {\em cut/cap\/} operations
depending on whether $w_1 w_2 w_3$ is two-sided or one-sided.
We call $T'$ a {\em pre-irreducible\/} triangulation.

If $w_1 w_2 w_3$ in $T$ is two-sided then cut along $w_1 w_2 w_3$ to produce
a surface $S_b$ with a boundary consisting of two disjoint 3-cycles 
$u'_1 u'_2 u'_3$
and $v'_1 v'_2 v'_3$ where $u'_i$ and $v'_i$ come 
from the original vertex $w_i$ for $i=1,2,3$. 
Cap the holes with two faces $u'_1 u'_2 u'_3$ and $v'_1 v'_2 v'_3$.
$T'$ is now a triangulation with $\mathrm{eg}(T') = \mathrm{eg}(T) - 2$.
If both $S_b$ and $S_f$ are orientable then $u'_1 u'_2 u'_3$ and 
$v'_1 v'_2 v'_3$ have opposite orientations in $S_b$.
If $S_b$ is orientable and $S_f$ is nonorientable then 
$u'_1 u'_2 u'_3$ and $v'_1 v'_2 v'_3$ have the same orientation in 
$S_b$.
If both $S_b$ and $S_f$ are nonorientable then $u'_1 u'_2 u'_3$ and 
$v'_1 v'_2 v'_3$ have no orientation in $S_b$.
It is not possible that $S_b$ is nonorientable and $S_f$ is 
orientable.

If $w_1 w_2 w_3$ in $T$ is one-sided then cut along $w_1 w_2 w_3$ to produce
a surface $S_b$ with a boundary consisting of the 6-cycle 
$u'_1 u'_2 u'_3 v'_1 v'_2 v'_3$ where $u'_i$ and $v'_i$
again come from the original vertex $w_i$ for $i=1,2,3$. 
Cap the hole with a new vertex $t'$ and six faces 
$t' u'_1 u'_2,t' u'_2 u'_3$, $t' u'_3 v'_1$,
$t' v'_1 v'_2$, $t' v'_2 v'_3$, and $t' v'_3 u'_1$.
$T'$ is now a triangulation with $\mathrm{eg}(T') = \mathrm{eg}(T) - 1$. 
Since $w_1 w_2 w_3$ is one-sided $S_f$ is nonorientable.
If $\mathrm{eg}(S_f)$ is odd then $S_b$ might be orientable or nonorientable.
If $\mathrm{eg}(S_f)$ is even then $S_b$ must be nonorientable.

Either cut/cap operation produces a pre-irreducible triangulation of 
$S_b$ with a smaller Euler genus which
by a sequence of edge contractions produces an irreducible triangulation with 
this smaller Euler genus.

The general plan of our algorithm is to reverse the procedure.
Start with the set of irreducible triangulations of $S_b$.
Generate other triangulations of $S_b$.
Check each generated triangulation to determine if it pre-irreducible for
a triangulation of $S_f$.
If it is then transform it into an irreducible triangulation of $S_f$ 
by applying the inverse of the cut/cap operation.

A simple algorithm for generating the pre-irreducible triangulations with $n$
vertices is

\begin{itemize}

\item $\mathcal{B}_s$ is the set of irreducible triangulations of $S_b$ 
with at most $n$ vertices;

\item $r_s$ is the generating rule ``if the triangulation has less than $n$ 
vertices split a vertex'';

\item $\overline{r_s}$ is ``if there is a contractible edge then contract a 
contractible edge'';

\item $\mathcal{T}_s$ is the set of triangulations of $S_b$ with at most 
$n$ vertices;

\item $f_s$ is the condition ``the triangulation is pre-irreducible'';

\item $\mathcal{F}_s$ is the set of triangulations of $S_b$ 
with at most $n$ vertices which are pre-irreducible.

\end{itemize}

By applying this algorithm for increasing values of $n$ all triangulations of 
$S_b$ which are pre-irreducible will be generated.
But from the algorithm there is no way of determining when all such 
triangulations have been generated.
Nakamoto and Ota \cite{MR1348564} provide a limit for the number of vertices
in the irreducible triangulations of $S_f$ but this limit does not provide
a practical criteria for termination.

This algorithm has been modified by making the generating rule more 
restrictive in two ways.
First, the order in which vertices are split is restricted.
Second, the generating rule is not applied to triangulations which are 
recognized as not being able to produce 
pre-irreducible triangulations. 

In the two-sided case the subgraph of $T'$ with the vertices and edges of the two
added faces is said to be the {\em handle frame\/}.
In the one-sided case the subgraph of $T'$ with the vertices and edges of the six
added faces is said to be the {\em crosscap frame\/}.
Vertices of the handle or crosscap frame are {\em frame vertices\/}.
The nonframe vertices are the {\em interior vertices\/}.
Edges of the handle or crosscap frame are {\em frame edges\/}.
Edges with two interior vertices are {\em interior edges\/}.
Edges with one interior vertex and one frame vertex are {\em support edges\/}.
Nonframe edges with two frame vertices are {\em crossframe edges\/} and
do not occur in pre-irreducible triangulations.

When reducing pre-irreducible triangulations to irreducible triangulations
by edge contraction the number of possible intermediate triangulations 
can be decreased by contracting edges in order: 
all contractible interior edges, all contractible 
support edges, and, finally, all contractible frame and crossframe edges.


As discussed above there may be as many as three different ways to generate
irreducible triangulations of $S_f$ depending the cut/cap operation used
and depending on the orientability of $S_b$.
For the presentation of the algorithm we assumed that $S_f$, $S_b$,
and the cut/cap operation are fixed.

The algorithm has three stages.
Using the reduction procedure as a guide we describe the stages in reverse
order.

\subsection{Stage 3}

Consider an irreducible triangulation $T$ and a pre-irreducible triangulation 
$T'$ obtained from $T$ using a cut/cap operation.
Let $xyz$ be a nonfacial 3-cycle in $T$.
The three edges in $xyz$ are each noncontractible.
Let $x'$ in $T'$ correspond to $x$ if $x \notin \{w_1,w_2,w_3\}$.
Let $y'$ and $z'$ be similarly defined.
There are several possibilities for what happens to $xyz$ in $T'$.

{\em Case 0:\/} Suppose $xyz$ and $w_1 w_2 w_3$ have no vertices in common.
Then the 3-cycle $x'y'z'$ in $T'$ is not a face and the three edges in 
$x'y'z'$ are each noncontractible.

{\em Case 1a:\/} Suppose $xyz$ and $w_1 w_2 w_3$ have exactly one vertex in 
common, say $x = w_1$,
and suppose $xy$ and $xz$ are on the same side of $w_1 w_2 w_3$, 
say $y'u'_1$ and $z'u'_1$ are edges of $T'$.
Then the 3-cycle $u'_1 y' z'$ in $T'$ is not a face and the three edges in 
$u'_1 y'z'$ are each noncontractible.

{\em Case 1b:\/} Suppose $xyz$ and $w_1 w_2 w_3$ have exactly one vertex in 
common, say $x = w_1$, 
and suppose $xy$ and $xz$ are on opposite sides of $w_1 w_2 w_3$, say 
$y'u'_1$ and 
$z'v'_1$ are edges of $T'$.
Then there is a path in $T'$ of length 3 from $u'_1$ to $v'_1$ containing the 
edges $u'_1 y'$, $y' z'$, $z' v'_1$.

{\em Case 2a:\/} Suppose $xyz$ and $w_1 w_2 w_3$ have exactly two vertices in 
common, say $x = w_1$ and $y = w_2$,
and suppose $xz$ and $yz$ are on the same side of $w_1 w_2 w_3$, 
say $z u'_1$ and $z u'_2$ are edges of $T'$.
Then the 3-cycle $u'_1 u'_2 z'$ in $T'$ is not a face and the three edges in 
$u'_1 u'_2 z'$ are each noncontractible.

{\em Case 2b:\/} Suppose $xyz$ and $w_1 w_2 w_3$ have exactly two vertices in 
common, say $x = w_1$ and $y = w_2$,
and suppose $xz$ and $yz$ are on opposite sides of $w_1 w_2 w_3$, 
say $z u'_1$ and $z v'_2$ are edges of $T'$.
Then there is a path in $T'$ of length 3 from $u'_1$ to $v'_1$ containing the 
edges $u'_1 z'$, $z' v'_2$, and $v'_2 v'_1$.
There is also a path in $T'$ of length 3 from $u'_2$ to $v'_2$ containing the 
edges $u'_2 u'_1$, $u'_1 z'$, and $z' v'_2$.

{\em Case 3:\/} Suppose $xyz$ and $w_1 w_2 w_3$ have three vertices in common. 
The corresponding edges in $T'$ are frame edges and
may or may not be contractible in $T'$.

Examining these cases we have:

For every edge of $T'$: 
\begin{itemize}
\item  the edge is a frame edge, or
\item  the edge is not contractible, or
\item  the edge is on a path of length 3 connecting $u'_i$ and $v'_i$ for 
some $i=1,2,3$.
\end{itemize}

Suppose $y'$ is an interior vertex of $T'$ which is not incident on any 
noncontractible edges.
Then by Theorem~\ref{transverse} the corresponding $y$ in $T$ has four ordered
neighbors $y_1$, $y_2$, $y_3$, and $y_4$ such that $yy_1y_3$ and $yy_2y_4$ are
nonfacial 3-cycles in $T$.
Either $y_1$ or $y_3$ must be on $w_1w_2w_3$ since $y'y'_1$ and $y'y'_3$ are
contractible, say $y_1 = w_i$. 
Likewise, either $y_2$ or $y_4$ must be on $w_1w_2w_3$, 
say $y_2 = w_j$.
Since $y_1 \neq y_2$ then $i \neq j$ and
$y'$ is on a path of length 3 connecting $u'_i$ and $v'_i$ and 
on a path of length 3 connecting $u'_j$ and $v'_j$ for $1 \leq i < j \leq 3$.

For every vertex of $T'$:
\begin{itemize}
\item  the vertex is a frame vertex, or
\item  the vertex is on a noncontractible edge, or
\item  the vertex is on a path of length 3 connecting $u'_i$ and $v'_i$ 
and on a path of length 3 connecting $u'_j$ and $v'_j$ for 
$1 \leq i < j \leq 3$.
\end{itemize}

As interior edges of $T'$ are contracted nonframe noncontractible edges
remain noncontractible, frame edges remain
frame edges, and paths connecting $u'_i$ and $v'_i$ do not become longer.

\begin{theorem}
\label{stageIII}
Let $T$ be an irreducible triangulation of a surface other than $S_0$.
Let $T'$ be a pre-irreducible triangulation 
obtained from $T$ using a cut/cap operation.
Let $T''$ be a triangulation obtained from $T'$ by contracting zero or more
interior edges.

For every edge of $T''$: 
\begin{itemize}
\item  the edge is a frame edge, or
\item  the edge is not contractible, or
\item  the edge is on a path of length at most 3 connecting $u'_i$ and 
$v'_i$ for some $i=1,2,3$.
\end{itemize}

For every vertex of $T''$:
\begin{itemize}
\item  the vertex is a frame vertex, or
\item  the vertex is on a noncontractible edge, or
\item  the vertex is on a path of length at most 3 connecting $u'_i$ 
and $v'_i$ 
and on a path of length at most 3 connecting $u'_j$ and $v'_j$ for 
$1 \leq i < j \leq 3$.
\end{itemize}
\end{theorem}

The stage 3 algorithm is

\begin{itemize}

\item $\mathcal{B}_3$ is the set of triangulations of $S_b$ 
which have a frame such that there are no contractible interior edges
and no crossframe edges;

\item $r_3$ is the generating rule ``split an interior vertex if after the 
split
every contractible nonframe edge is on a path of length at most 3 connecting 
$u'_i$ and $v'_i$ for some $i=1,2,3$ and 
every nonframe vertex that is not on a noncontractible edge
is on a path of length at most 3 connecting $u'_i$ 
and $v'_i$ 
and on a path of length at most 3 connecting $u'_j$ and $v'_j$ for 
$1 \leq i < j \leq 3$'';

\item $\overline{r_3}$ is ``if there is a contractible interior edge then 
contract a contractible interior edge'';

\item $\mathcal{T}_3$ is the set of 
all triangulations for which 
every contractible nonframe edge is on a path of length at most 3 connecting 
$u'_i$ and $v'_i$ for some $i=1,2,3$ and 
every nonframe vertex that is not on a noncontractible edge
is on a path of length at most 3 connecting $u'_i$ 
and $v'_i$ 
and on a path of length at most 3 connecting $u'_j$ and $v'_j$ for 
$1 \leq i < j \leq 3$;

\item $f_3$ is the condition ``every path connecting 
$u'_i$ and $v'_i$ for $i=1,2,3$ has length at least 3'';

\item $\mathcal{F}_3$ is the set of all pre-irreducible triangulations of 
$S_b$ which produce irreducible triangulations of $S_f$ using the
inverse of the cut/cap operation.

\end{itemize}

\begin{question}
Is $\mathcal{T}_3$ finite for all $S_b$?
\end{question}

\subsection{Stage 2}

In stage 2 we want to generate $\mathcal{F}_2 = \mathcal{B}_3$. 
This is done by splitting frame vertices to produce support edges
without altering the frame.

The stage 2 algorithm is

\begin{itemize}

\item $\mathcal{B}_2$ is the set of triangulations of $S_b$ 
which have a frame such that there are no contractible interior edges
and no contractible support edges;

\item $r_2$ is the generating rule ``split a frame vertex if after the 
split the frame is unchanged and there are no contractible interior edges'';

\item $\overline{r_2}$ is ``if there is a contractible support edge then 
contract a contractible support edge'';

\item $\mathcal{T}_2$ is the set of triangulations of $S_b$ 
which have a frame such that 
every contractible nonframe edge is a support edge or a crossframe edge;

\item $f_2$ is the condition ``there are no crossframe edges'';

\item $\mathcal{F}_2$ is the set of triangulations of $S_b$ 
which have a frame such that there are no contractible interior edges
and no crossframe edges.

\end{itemize}

\begin{question}
Is $\mathcal{T}_2$ finite for all $S_b$? 
\end{question}

\subsection{Stage 1}

In stage 1 we want to generate $\mathcal{F}_1 = \mathcal{B}_2$. 
This is done by splitting any vertex to produce the frame.
Since the two types of frames are different there are two stage 1 algorithms.

The stage 1 algorithm of two-sided cut is

\begin{itemize}

\item $\mathcal{B}_1$ is the set of irreducible triangulations of $S_b$;

\item $r_1$ is the generating rule ``split a vertex if after the 
split there are at most six contractible vertices'';

\item $\overline{r_1}$ is ``if there is a contractible edge then contract a 
contractible edge'';

\item $\mathcal{T}_1$ is the set of triangulations of $S_b$ 
which have at most six contractible vertices;

\item $f_1$ is the condition ``there are two faces with no common vertices
which contain all the contractible vertices'';

\item $\mathcal{F}_1$ is the set of triangulations of $S_b$ 
which have a frame such that there are no contractible interior edges
and no contractible support edges.

\end{itemize}

The stage 1 algorithm of one-sided cut is

\begin{itemize}

\item $\mathcal{B}_1$ is the set of irreducible triangulations of $S_b$;

\item $r_1$ is the generating rule ``split a vertex if after the 
split there are at most seven contractible vertices and there is a vertex
$t$ such that the set of $t$ and its neighbors contain all of the contractible
vertices'';

\item $\overline{r_1}$ is ``if there is a contractible edge then contract a 
contractible edge'';

\item $\mathcal{T}_1$ is the set of triangulations of $S_b$ 
which have at most seven contractible vertices and have a vertex
$t$ such that the set of $t$ and its neighbors contain all of the contractible
vertices;

\item $f_1$ is the condition ``there is a vertex $t$ of degree six 
such that the set of $t$ and its neighbors contain all of the contractible
vertices'';

\item $\mathcal{F}_1$ is the set of triangulations of $S_b$ 
which have a frame such that there are no contractible interior edges
and no contractible support edges.

\end{itemize}

Each application of $r_1$ increases the number of contractible vertices. 
There can be at most seven contractible vertices,
therefore, $\mathcal{T}_1$ is finite and $\mathcal{F}_1$ is finite.

This three stage algorithm was implemented as two computer programs 
{\em sgh} (Surftri Grow Handle) and {\em sgc} (Surftri Grow Crosscap).

{\em sgh} and {\em sgc} terminate when $S_f$ is $S_1$, $S_2$,
$N_1$, $N_2$, $N_3$, or $N_4$ which means for these surfaces $\mathcal{T}_2$ 
and $\mathcal{T}_3$ are finite.
The irreducible triangulations of these surfaces are available as computer 
files \cite{surftri}.
The approximate total computer times required on a cluster of computers with 
an average CPU speed of about 2 GHz were: 
under 1 second for $S_1$, $N_1$, and $N_2$;
30 minutes for $N_3$ (5.4 irreducible triangulations per second);
34 hours for $S_2$ (3.2 irreducible triangulations per second); and
56 days for $N_4$ (1.3 irreducible triangulations per second).
One is tempted to extrapolate for $S_3$ and $N_5$ and obtain times in 
centuries.

If there is a general formula for a limit on the number of vertex 
splittings required to produce pre-irreducible
triangulations of $S_f$ from irreducible triangulations of 
$S_b$  then the algorithm terminates for all surfaces.
This formula might also lead to a smaller upper limit on the number of vertices
in the irreducible triangulations of $S_f$ than the limit in 
\cite{MR1348564}.

\section{Verification}
\label{verify}

\subsection{Known results}

The irreducible triangulations of $S_1$, $N_1$, and $N_2$ have previously been
determined.  
{\em sgh} and {\em sgc} produces these same triangulations.

For each of the surfaces for which we have irreducible triangulations 
{\em surftri\/} has been used to generate and count all of the triangulations 
up to some number of vertices.
These counts are available on line \cite{surftri}.

Counts of the number of triangulations of $S_0$ are available from 
{\em plantri\/} \cite{plantri}.  
These same counts were produced by {\em surftri\/} for
triangulations up to 15 vertices and for triangulations with minimum degree 4
up to 18 vertices.
{\em surftri\/} runs at about half the rate of {\em plantri\/}.
In order to handle nonorientable triangulations the data structure for
triangulations in {\em surftri\/} is twice the size as in {\em plantri\/}.
The rate difference is due to the additional processing required 
to manipulate this larger data structure and 
to maintain information about which edges are contractible.

Counts of the number of triangulations of $S_1$ are available from 
{\em togen\/} \cite{MR1912863}.
Agreement was found for up to 14 vertices.

Lutz \cite{math.CO/0506316} \cite{math.CO/0506372} has generated all the 
triangulations of surfaces
up to 10 vertices by starting with one face and adding one face at a time in
lexicographic order.
This method has been extended by Sulanke and Lutz \cite{isofree}
to generate triangulations of surfaces up to 12 vertices.  
The counts of triangulations of surfaces for which we have irreducible
triangulations agree.

\subsection{Redundancy}

{\em sgh} and {\em sgc} together generate all the pre-irreducible 
triangulations of $S_f$.  
For example, when $S_f = S_2$ the average number of pre-irreducible 
triangulations per irreducible triangulation is about 33.
We take advantage of this redundancy.
We count the nonseparating 3-cycles in each irreducible triangulation generated
and check that this number agrees with the number of pre-irreducible 
triangulations generated.

\subsection{Random search}

A method for searching for irreducible triangulations of $S_2$ is to generate
random triangulations and randomly contract edges until no more edges can be
contracted.  The starting random triangulations used had 20 vertices. 
Random irreducible triangulations were generated at the rate of about 1,000 
per second for six days.  About 97\% of the irreducible triangulations
generated by {\em sgh} were found.
This method assumes that all irreducible triangulations have 20 or fewer 
vertices and that $N(S_2) <= 20$ where $N(S)$ is defined in the next 
subsection.

A similar random search found 9705 of the 9708 irreducible triangulations 
of $N_3$.

\subsection{Start with pseudo-minimal}

Two triangulations $T$ and $T'$ of
a surface are {\em equivalent\/} 
if there is a isomorphism $h$
with $h(T)=T'$.
That is, if $a$, $b$, and $c$ are vertices of $T$ then
$ab$ is an edge of $T$ if and only if $h(a)h(b)$ is an edge
of $T'$ and a face of $T$ is bounded by the cycle $abc$
if and only if a face of $T'$ is bounded by the cycle
$h(a)h(b)h(c)$.

Let $ac$ be an edge in a triangulation $T$ and
$abc$ and $acd$ be the two faces which have $ac$ as a common edge.
The {\em diagonal flip\/} of $ac$ is obtained by deleting $ac$,
adding edge $bd$, deleting the faces $abc$ and $acd$, and adding the faces
$abd$ and $bcd$.
Two triangulations are {\em equivalent under diagonal flips\/}
if one is equivalent to a triangulation obtained from the other 
by a sequence of diagonal flips.

The number of vertices of an irreducible triangulation 
can not be reduced by edge contraction.
Negami \cite{MR95m:05091} defines a type of triangulation for which 
the number of vertices can not be reduced by a combination of diagonal flips 
and edge contractions. 
An irreducible triangulation is said to be {\em pseudo-minimal\/} 
if it is equivalent under diagonal flips only to irreducible triangulations.

Define $N(S)$ to be the minimum value such that two triangulations $T$ and $T'$
of $S$ are equivalent under diagonal flips 
if the number of vertices in $T$ and the number of vertices in $T'$ 
are equal and 
at least $N(S)$.  Negami \cite{MR95m:05091} has
shown that such a finite value exists for any $S$.

We can start with the 865 pseudo-minimal triangulations of $S_2$ and use
vertex splitting and diagonal flips to generate irreducible triangulations
of $S_2$.
When using this method we assume that we know all of the pseudo-minimal 
triangulations and that the maximum number of vertices
in any irreducible triangulation is 17.  This value of 17 is the number of
vertices in the irreducible triangulation obtained by joining one face of
each of two copies of the
largest irreducible triangulation of $S_1$.
Approximately $10^{13}$ triangulations with minimum degree of 4 were generated
using about 5 years of computer time.
This method demonstrates only the self consistency of the irreducible 
triangulations of $S_2$ since we use the set of irreducible 
triangulations to find the set of pseudo-minimal 
triangulations and we use the set of pseudo-minimal 
triangulations to find the set of irreducible 
triangulations.

This same method was also used to generate all the irreducible triangulations 
of $N_3$ starting with the 133 pseudo-minimal triangulations and generating 
all ($> 10^{12}$) triangulations with at most 17 vertices and minimum 
degree 4.

\subsection{Generating triangulations by diagonal flip}

Let $n \geq N(S)$ then we can construct a triangulation of $S$ 
which has $n$ vertices, e.g. using the random method described above.
From this triangulation we generate by diagonal flips all the triangulations 
of $S$ with 
$n$ vertices using backtracking and recording all the generated triangulations.
We compare these triangulations with those generated by {\em surftri\/}.
For every irreducible triangulation of $S$ with at most $n$ vertices
there are triangulations with $n$ vertices which can be obtained by vertex 
splitting only from this one irreducible triangulation.
This method checks both {\em surftri\/} and the irreducible triangulations.

$N(S_0)$, $N(S_1)$, $N(N_1)$, and $N(N_2)$ are known \cite{wagner} 
\cite{MR46:8878} \cite{MR91g:05038}.
Checking the irreducible triangulations generated for $S_2$, $N_3$, and $N_4$
we have determined \cite{properties} that
$N(S_2) = 10$, $N(N_3) = 9$, and $N(N_4) = 10$ thus this method demonstrates 
only the self consistency of the irreducible triangulations of these three
surfaces.
For performance reasons we stored the triangulations in memory and restricted
this storage to 7 gigabytes.
Agreement was found up to 16 vertices for $S_0$, 13 for $S_1$, 12 for $S_2$, 
14 for $N_1$, 13 for $N_2$, 12 for $N_3$, and 11 for $N_4$.

A similar method, restricting our consideration to those triangulations with
minimum degree of 4, allows for triangulations with more vertices.
Agreement was found up to 19 vertices for $S_0$, 14 for $S_1$, 12 for $S_2$, 
16 for $N_1$, 14 for $N_2$, 12 for $N_3$, and 11 for $N_4$.


\bibliographystyle{amsplain}
\bibliography{gen_triang}

\end{document}